\documentclass[a4paper,12pt,reqno]{amsart}

\usepackage[applemac]{inputenc}
\usepackage[T1]{fontenc}
\usepackage{graphicx}

\usepackage{amssymb, microtype,graphicx}
\usepackage{booktabs}  

\newcounter{notecounter}

\newtheorem{Th}{Theorem}

\theoremstyle{definition}

\title{Using Hoare's Theorem to find the signature of a subgroup of an NEC group.}
\author{ David Singerman and Paul Watson}
\date{\today}

\begin{document}
\maketitle

\section{Introduction}

There is a standard method, due to the first author, \cite{S1} of finding the signature of a subgroup $\Lambda$ of a Fuchsian group $\Gamma$ given the permutationn representation of $\Gamma$ on the $\Lambda$-cosets. The structure of a non-Euclidean crystallographic (NEC) group is considerably more complicated, but a method of computing the signature of a subgroup was found by A.~H.~M.~ Hoare  in \cite{H} The purpose of this article is to explain how Hoare's Theorem is used and to give examples. Partial techniques to find the signature of a subgroup a are due to Gromadzki, \cite{G}.

The signature of a cocompact NEC group $\Gamma$ has the form 

$$(g;\pm; [m_1,\ldots m_r];\{(n_{11},\ldots n_{1s_1}),\ldots \{(n_{k1},\ldots n_{ks_k})\}) \eqno(1)$$

Here $g$ is the genus of the quotient space $\mathbb{H}/\Gamma$, $(\mathbb{H}$ being the hyperbolic plane), $m_1,\ldots,m_r$ are the proper periods and the $n_{ij}$ are the {\it link} periods;  the brackets ($n_{u_1}\ldots n_{u_{s_u}})$ are called the {\it period cycles.}

Our problem is to find the signature of a subgroup $\Lambda$ of finite index in $\Gamma$, given the permutation representation of $\Gamma$ on the $\Lambda$ cosets.

If $\Gamma$ is a Fuchsian group then the result is in \cite{S1}. For $\Gamma$ a proper NEC group the result follows from Hoare's theorem in \cite{H}. This is rather more complicated to apply than for Fuchsian groups and the purpose of this note is to give a fairly algorithmic procedure to apply Hoare's result. (Much of this is derived from results in the second author's Southampton Ph.D. thesis,\cite{W}.)

The proper periods of the subgroup are computed in the same way as for Fuchsian groups. The major problem is to compute the link periods and the period cycles for $\Lambda$.

All reflections of $\Lambda$ are conjugate to one of the generating reflections of $\Gamma$ and so these correspond to fixed points of the generating reflections of $\Gamma$ when acting on the right cosets of the subgroup by right multiplication. Thus if a reflection generator $c_i$ of $\Gamma$ fixes a coset $\Lambda g_j$ then $c_{ij}=g_jc_ig_j^{-1} \in \Lambda.$ We call $c_{ij}$ an {\it induced reflection} of $\Lambda.$ Now suppose that we have a period cycle $(n_{i1},n_{i2}\ldots n_{is_i})$.  This corresponds to a part of the presentation of $\Gamma$ which has generators $c_{i0},\ldots c_{is_i}\,, e$ obeying the relations 

$$c_{i{j-1}}^2=c_{ij}^2=(c_{i{j-1}}c_{ij})^{n_{ij}}=1= e_ic_0e_i^{-1}c_s =1  \eqno(2)$$

We also have elliptic generators $x_i$ of orders $m_i$ so we have relations

$$x_j^{m_j}=1  \eqno(3)$$ and also the long relations  

$$x_1\ldots x_r e_1 \ldots e_k a_1b_1a_1^{-1}b_1^{-1}\ldots a_gb_ga_g^{-1}b_g^{-1}=1     \eqno(4) $$     
if $\mathbb{H}/\Gamma$ is orientable or 
$$x_1\ldots x_r e_1 \ldots e_k  a_1^2\ldots a_g^2=1     \eqno(5)    $$              
if  $\mathbb{H}/\Gamma$ is non-orientable.

The group $\Gamma$ has a presentation which consists of the generators $c_1,\ldots, c_k,e_1,\ldots ,e_k, a_1,b_1,\ldots ,a_g,b_g$

if $\mathbb{H}/\Gamma$ is orientable  and $c_1\ldots ,c_k,e_1, \ldots,e_k, a_1,\ldots ,a_g$ 

if $\mathbb{H}/\Gamma$ is non-orientable, and relations (2), (3),(4) if $\mathbb{H}/\Gamma$ is orientable and (2),(3), (5) if $\mathbb{H}/\Gamma$ is non-orientable. 

The elements $c_i$ are reflections, $x_i$ elliptics, $a_i,b_i$ hyperbolic in the orientable case and the $a_i$ are glide-reflections in the non-orientable case. The $e_i$, called the connecting generators, are orientation-preserving and nearly always hyperbolic.

For these presentations see \cite{M}, \cite {BEGG}.

{\bf Definition.} Two reflections $c$ and $d$ of $\Gamma$ are called {\it linked}
with {\it link period $n$ } if $cd$ is an elliptic element of order $n$.

The following ideas come from Hoare's paper \cite{H}.

If $c$ and $d$ are linked reflection generators with link period $n$ then $c$ and $d$ generate a dihedral group $D_n$ of order $2n$. Let $\sigma$ be an orbit of the $\Lambda$-cosets under $D_n$ and let $K$ be a coset in $\sigma$. If $m$ is the least positive integer such that $(cd)^m \in K$ then either;

a) $\sigma$ contains no coset fixed by $c$ or $d$ in which case $\sigma$ has length $2m$ and gives an elliptic generator for $\Lambda$ with period $n/m$, or

b) $\sigma$ contains two cosets fixed by $c$ and $d$, one fixed by each if $m$ is odd, two fixed by one and none by the other if $m$ is even. Now $\sigma$ has length $m$ and the refection generators  of $\Lambda$ corresponding to the two fixed cosets are linked with period $n/m$. Each reflection generator of $\Lambda$ appears in precisely two of these links, unless it is linked only to itself. This gives the period cycles of $\Lambda$ each coming from one of the period cycles of $\Gamma.$

c)  In a) above we find out how to get some proper periods of $\Lambda$. Other proper periods of $\Lambda$ are just those induced by the proper periods of $\Gamma$. These are found in exactly the same way as in \cite{S1}. That is if $x$ is an elliptic element of $\Gamma$ of period $n$, then we let $x$ act on the right $\Lambda$-cosets. by right multiplication. If we have an orbit (i.e. a cycle) of length $m$ then there is an induced elliptic period equal to $n/m$ in $\Lambda$. All elliptic periods of $\Lambda$ are found in a) or c).

d)  The orientability of $\Lambda$ is found by using \cite{HS}.  If $\Gamma$ is a group with generators $\Phi$ and $\Lambda$ is a subgroup then the Schreier coset graph ${\mathcal H}(\Gamma,\Lambda,\Phi)$  is the graph with vertices the cosets of $\Lambda$ in $\Gamma$ and directed edges at each vertex for each $b\in \Phi$ such that $b:\Lambda g\longrightarrow \Lambda gb$. If $c$ is a reflection and if $\Lambda gc=\Lambda g$ then the directed edge $c:\Lambda g\longrightarrow \Lambda gc$ is a {\it reflection loop}. Let $\overline{\mathcal H}$ be the Schreier graph with the reflection loops deleted. We call this the {\it augmented Schreier graph.} Each path in  $\overline{\mathcal H}$  corresponds to a word in $\Phi$ and hence to an element of $\Gamma$. A path is called {\it positive (negative)} if it corresponds to a orientation-preserving (orientation- reversing ) element of $\Gamma$. We label the edges of this graph with the generators of the group.

\bigskip

Then it was shown in \cite{HS} that $\mathbb{H}/\Lambda$ is orientable if and only if all circuits in  $\overline{\mathcal H}$ are positive, where $\mathbb{H}$ is the hyperbolic plane.

\bigskip

Finally, the genus of $\mathbb{H}/\Lambda$ is calculated using the Riemann-Hurwitz formula as usual.  Usually, it is best to find the index $|\Gamma^+:\Lambda^+|$ where $\Gamma^+$ (resp. $\Lambda^+)$ is the canonical Fuchsian group of $\Gamma$, (resp. $\Lambda$) the index two subgroups of orientation-preserving isometries.  The signature of $\Gamma^+$ is found from \cite{S2}.

\section {The algorithm.}

We suppose that $\Gamma$ has a subgroup of index $N$ and let the cosets be 
$\Lambda g_1, \Lambda g_2,\ldots\Lambda g_N$, which we represent in the usual way by $1,2,\ldots.N$.


Then 

{\it Step 1.}  Write down the induced reflections.

\medskip

{\it Step 2.}  Find all the links and link periods. Typically, $c=c_{i{j-1}}$and $d=c_{ij} $ are linked with link period $n=n_{ij}$.  Find the image of $cd$ in the symmetric group $S_N$ and note the cycle lengths in this permutation. In particular the the lengths of the cycles containing fixed points of $c$ or $d$ will give the integers $m$ above.

\medskip

{\it Step 3.} Let $c$ and $d$ be two linked refections with link period $n$.  Then $c$ and $d$ generate a dihedral group $D_n$. Find the orbits of this $D_n$ in its action on the $\Lambda$-cosets. Let $\sigma$ be an orbit and let $K$ be a coset in $\sigma$ (So $K$ is just one of the integers $1, 2, \cdots ,N$.) Consider the cycles of $cd$. If there is a cycle of length $m$ this just means that $K(cd)^m=K$. 

\medskip 

{\it step 4}. Find the periods and link periods of the subgroup. If there are no cosets fixed by $c$ or $d$ in $\sigma$ then Hoare's results in  \cite{H} tell us that $\sigma$ has length $2m$ and gives an elliptic period $n/m$ in $\Lambda$.

If there are fixed cosets then 

 $\sigma$ contains exactly two cosets fixed by $c$ and $d$, one fixed by each if $m$ is odd, two fixed by one and none by there other if $m$ is even. Now $\sigma$ has length $m$ and the induced reflection generators of $\Lambda$ are linked  with link period $n/m$.  We write $c$  linked to $d$ by $c\sim d$. Each reflection generator of $\Lambda$ occurs in exactly two of these links  unless it is linked only to itself. This gives the period cycles of $\Lambda$, each coming from one of the period cycles of $\Gamma$. 
 
 \medskip
 
 {Step 5.} How to deal with the relation $ec_0e^{-1}c_s=1$.
 
 Regard $\langle{ec_0e^{-1}, c_s\rangle}\cong D_1$. 
 Now the orbits of $\langle{ec_0e^{-1}, c_s}\rangle\cong D_1$ are the two-cycles of $c_s=ec_0e^{-1}$ or the one-cycles, (fixed points). If we have a 2-cycle $(\gamma,\delta)$ then $\{\gamma,\delta\}$ is an orbit containing no fixed points of $c_s$ or 
 $ec_0e ^{-1}(=c_s)$. Now $ec_0e^{-1}c_s=1$ and so $m=1$ and there is an
elliptic period equal to $1/1=1$, that is no elliptic period at all. 
If $c_s=ec_0e^{-1}$ fixes $k$ then $k(ec_0e^{-1})=k$ so $c_0$ fixes $ke$ and we have a link $c_s\sim c_{0ke}$ with link period 1.

{\bf Example 1}
 
In our first examples we find all subgroups of index two in an extended triangle group. This is a group $\Delta$ generated by three reflections in the sides of a triangle with angles $\pi/n_1,\pi/n_2, \pi/n_3$.  This has signature $(0;+;[ \ ];\{n_1,n_2,n_3\})$, which we usually denote by $(n_1,n_2,n_3)$,  which has presentation 
 
$ \{c_1,c_2,c_3|c_1^2=c_2^2=c_3^2=(c_1c_2)^{n_1}=(c_2c_3)^{n_2}=(c_3c_1)^{n_3}=1\}.$

\bigskip

There  are 7 epimorphisms $\theta_1:\Delta\longrightarrow C_2=\{e,t\}$, ($t^2=1$).

\begin{enumerate} 

\item $\theta_1(c_i)=t$ for $i=1,2,3$. Here the kernel (equal to the stabiliser of a point) is the canonical Fuchsian triangle group with signature  $(0;+;[n_1,n_2,n_3];\{\ \})$ usually denoted by $[n_1,n_2,n_3]$.

\medskip

\item $\theta_2(c_1)=\theta(c_2)=t=(1,2)$, $\theta (c_3)=e=(1)(2)$.
Note that as $\theta(c_1c_3)=\theta(c_2c_3)=t$, $n_2$ and $n_3$ must be even.
\end{enumerate}
\medskip

We now use our algorithm to find $\Lambda$ the stabiliser of a point. 

\medskip

Step 1. The induced reflections are $c_{31}, c_{32}$. 

\medskip

Step 2. Find the links and link periods. Here the links are $c_1\sim c_2$  with link period $n_1$, $c_2\sim c_3$ with link period $n_2$ and $c_3\sim c_1$ with link period $n_3$.

\medskip

Step 3. Find the orbits of the dihedral subgroups. There are three dihedral subgroups here $\langle{c_1,c_2\rangle}\cong D_l$,  $\langle{c_2,c_3}\rangle \cong D_m$, and $\langle{c_3,c_1}\rangle \cong D_n$. The orbits of all these three are $\{1,2\}$.  Also $c_1c_2=(1)(2)$, and $c_2c_3=c_3c_1=(1,2)$.

Step 4. We first notice that no coset is fixed by $c_1, c_2$. The orbit $\sigma$ of 
$\langle{c_1,c_2\rangle}$ has length 2 so that here $m=1$ and so there is an elliptic period $n_1$ in $\Lambda$.

Now we consider the dihedral group $\langle{c_2,c_3}\rangle$.  Here the orbit has length 2 (which is even) and $c_3$ fixes two cosets 1 and 2 and $c_2,$ fixes no coset so $c_{31}\sim c_{32}$. As $c_1c_3=(1,2)$, $m=2$ and so the link period is $n_2/2$. Similarly , we also get  $c_{31}\sim c_{32}$ from the orbit $\langle{c_3,c_1}\rangle$. Here the link period is $n_3/2$. We get one chain $c_3\sim c_2\sim c_3$ and the period cycle is $\{n_2/2,n_3/2\}$. As the augmented Schreier graph has two vertices, it cannot have any negative paths and so $\Lambda$ has orientable quotient space.  Thus the signature of $\Lambda$ is of the form $(g;+; [n_1];(\{n_2/2,n_3/2\})$. We compute $g$ using  Riemann-Hurwitz. We pass to the canonical Fuchsian groups.   $\Delta^{+}$ has Fuchsian signature $(0; (n_1,n_2,n_3))$ and $\Lambda^{+}$ has Fuchsian signature $(2g+1-1;n_1,n_1,n_2/2,n_3/2)$. As the index of $\lambda^+$ in $\Gamma^+$ is two we apply Riemann-Hurwitz to get

$$4g-2+1-1/n_1+1-1/n_1+1-2/n_2+1-2/n_3=2(1-1/n_1-1/n_2-1/n_3)$$

giving $g=0$. Thus the signature of the NEC group $\Lambda$ is

$$(0; +; [n_1]; \{(n_2/2,n_3/2)\}). \eqno(6)$$

$\theta_2(c_1)=(1,2)$,  $\theta_2(c_2)=\theta(c_3)=(1)(2)$.
Again we let $\Lambda$ denote the stabiliser of a point. 

Step 1. The induced reflections are $c_{21}, c_{22}, c_{31}, c_{32}$

Step 2. The links are as in the previous example. 
 
 Step 3. The orbits of $\langle{c_1,c_2}\rangle$  and $\langle{c_1, c_3}\rangle$ are $\{1,2\}$.  The orbits of $\langle {c_2,c_3}\rangle$ are $\{1\}$ and $\{2\}$.  Also $c_1c_2=(1,2)$
 
 Step 4. The orbit of $\langle{c_1,c_2}\rangle$   contains  two fixed points of $c_2$ and none of $c_1$.  As $c_1c_2=(1,2)$ we find a link $c_{21}\sim c_{22}$ with link period $n_1/2$.
 
 Now $c_2c_3=(1)(2)$. By considering the two singleton orbits of $\langle {c_2,c_3}\rangle$ , we find that the orbit $\{1\}$ contains 1 fixed point of $c_2$ and one fixed point of $c_3$. Here $c_2c_3$ is the identity so $m=1$ (odd!) so we get a link $c_{21}\sim c_{31}$ with link period $n_2$. Similarly, by considering the orbit $\{2\}$ we find a link $c_{22}\sim c_{32}$ with link period $n_2$.
 
 Also, the orbit of $\langle c_1,c_3\rangle$ contains two fixed points of $c_3$ and none of $c_1$. (Note $c_1c_3=(1,2)$ so $m=2$, (even!)).  Thus $c_{31}\sim c_{32}$ with link period $n_3/2$. We thus have one chain
 
 $$c_{21}\sim c_{22}\sim c_{32}\sim c_{21}\sim c_{31}\sim c_{21}$$
 with link periods 
 
 $$n_1/2,n_2,n_3/2,n_2.$$
 
 Thus $\Lambda$ contains one period cycle $(n_1/2,n_2,n_3/2,n_2)$.
 
 As in the previous example, we show that  we show that $\Lambda$ has genus zero and orientable quotient space and so $\Lambda$ has signature 
 
 $$(0;+;\{(n_1/2,n_2,n_3/2,n_2\}). \eqno(7)$$
 
 By permuting $n_1,n_2, n_3$ we find the other possible subgroups of index two  leading to the extended triangle group  $(0;+;[\  ];\{n_1,n_2,n_3)\})$, to give 
 \begin{Th} The NEC triangle group $\Delta$  of signature $$(0;+;[\  ];\{n_1,n_2,n_3)\})$$  has the Fuchsian triangle group $[n_1,n_2,n_3]$ as a subgroup of index two. If $n_1,n_3$ are even then the  only other possible subgroups of index two in $(0;+;[\ ];\{n_1,n_2,n_3)\}$, have signatures
 
 \item $(0;+;[n_2];\{(n_2/2,n_3/2)\}$ 
\item $(0;+;;[\ ];\{(n_1/2,n_2,n_3/2,n_2\})$.

Similarly, if $n_1, n_2$ are even we get two other NEC subgroups of index 2 and if $n_2,n_3$ are even we get another two. Altogether this gives us our seven subgroups of index two in the extended triangle group, $(n_1,n_2,n_3)$.
\end{Th}

\bigskip

We now give other examples which illustrates all of the points used in Hoare's Theorem.  The first one  illustrates more fully how to find the period cycles and deal  with empty period cycles.

\bigskip

{\bf Example 2} 

Let $\Gamma$ be an NEC group of signature

 $$(0;+;[ \ ], \{(2,3),(\ )\}) \eqno(8)$$

The canonical presentation of $\Gamma$ is 

$$\langle c_0,c_1,c_2,d, e_1,e_2|c_0^2=c_1^2=(c_0c_1)^2=(c_1c_2)^3=d^2=1,
e_1c_0e_1^{-1}=c_2,e_2de_2^{-1}=d,e_1e_2=1\rangle.$$

\rm Consider the following permutation representation.

\medskip

$c_0\longmapsto (1,2)(3)(4)$

$c_1\longmapsto (1)(2)(3,4)$

$c_2\longmapsto (1,3)(2)(4)$

$e_1\longmapsto (1)(4)(2,3)$

$e_2\longmapsto (1)(4)(2,3)$

$d\longmapsto (1)(4)(2,3)$

\bigskip

We now follow our algorithm. 

Step 1. The induced reflections are $c_{03}, c_{04}, c_{11}, c_{12},  c_{22}, c_{24},d_{11},
d_{14}.$

\bigskip

Step 2. The links are $c_0\sim c_1$ with link period 2, $c_1\sim c_2$ with link period 3,
$e_1c_0e_1^{-1}\sim c_2$ with link period $1$, $e_2de_2^{-1}\sim d$ with link period 1.

\bigskip

Step 3. Find the orbits of the dihedral subgroups.  These are 
$\langle c_0,c_1\rangle \cong D_2.$ The orbits here are $\{1,2\}$,$\{3,4\}$, 
and $c_0c_1=(1,2)(3,4).$

\medskip

$\langle c_1,c_2\rangle \cong D_3$ The orbits are $\{1,3,4\}$, $\{2\}$ and $c_1c_2=(1,3,4)(2).$

\medskip

$\langle e_1c_0e_1^{-1},c_2\rangle \cong D_1$,  The orbits are $\{1,3\}, \{2\}$,$\{4\}$ and $ec_0e^{-1}c_2=(1)(2)(3)(4)$.

\medskip

  $\langle e_2de_2^{-1},d \rangle \cong D_1$. The orbits are $\{1,3\}, \{2\}$,$\{4\}$ and $ede^{-1}d=(1)(2)(3)(4)$.
  
  \bigskip
  
  Step 4. We have $c_0c_1=(1,2)(3,4)$. In the cycle $(1,2)$ we have two fixed points of $c_1$ and none of $c_0$. (Note $m=2$ here is even so we should have two fixed points of one generator and none of the other.)  As $n=2$ we have a link 
  $c_{11}\sim c_{12}$ with link period  $2/2=1$
  
  \medskip
  
  On the orbit $\{3,4\}$ of $\langle c_0,c_1\rangle$ we similarly get a link $c_{03}\sim c_{04}$ with link period equal to 1.
  
  \medskip
  
 Now $c_1c_2=(1,3,4)(2)$. In the cycle $(1,3,4)$ we have one fixed point 1 of $c_1$ and one fixed point 4 of $c_2$ so we get a link $c_{11}\sim c_{24}$  with link period $3/3=1$.
  
  \medskip
  
  On the cycle $\{2\}$ of $\langle c_1,c_2\rangle$ we have one fixed point of $c_1$ and one of $c_2$, namely 2 in both cases and so we get a link $c_{12}\sim c_{22}$ with link period $3/1=3$. 
  
  \medskip
  
  Step 5. We now consider the dihedral group  $\langle e_1c_0e_1^{-1}, c_2\rangle\cong D_1$ of order two.   We have two fixed points of $c_0,$ namely 2 and 4. From step 5, we get links $c_{22}\sim c_{02e_1}$ or $c_{22}\sim c_{03}$ with link period 1. Similarly we have the link $c_{24}\sim c_{04}$ with link period 1.
  
  \bigskip
  
  Putting all these links together we find that we only have one chain:
  $$c_{11}\sim c_{12}\sim c_{22}\sim c_{03}\sim c_{04}\sim c_{24}\sim c_{11}$$
  
  All the link periods are equal to one except for the link period between $c_{24}$ and $c_{14}$ which is equal to 3. We can omit the link periods equal to one and so we are left with a single period cycle $\{(3)\}$.
  
  We now find the number of period cycles induced by the empty period cycle of $\Gamma$. This empty period cycle corresponds to the reflection $d$, which gives two induced reflections $d_1$ and $d_4$. This gives us two empty period cycles of $\Lambda$
  
  We now consider the orientability of $\Lambda$.  So we consider the augmented coset graph as described above. This graph has 4 vertices 1,2,3,4.  We have an edge from 1 to 2, because of the  reflection $c_0$, an edge from 1 to 2 coming from the reflection $c_1$ and an edge from 1 to 3 because of the reflection $d$. Thus we have a triangle $1\longrightarrow 2\longrightarrow3 \longrightarrow 1.$
  This gives an orientation-reversing loop and so $\Lambda$ has non-orientable quotient space. 
  
  We now compute the genus of $\Lambda$. The group $\Gamma$ has signature (8).  $\Lambda$ has a signature of the form $(g;-[\  ];\{(3), (\ ),(\ ))$.  By [\cite{S2}$,\Gamma^+$  has Fuchsian signature $(1;2,3)$ and $\Lambda^+$ has signature  (g+2; 3).
  Riemann-Hurwitz now gives $2(g+2)-2+2/3=4.(7/6)=14/3$ and thus $g=1$ and so the signature of $\Lambda$ is $(1;[\ ];\{(3),(\ ),(\ )\})$.
  
 The next example which comes from \cite{W} illustrates more fully how to determine the proper periods of the subgroup. 
  
 \bigskip
 
 {\bf Example 3.} 
 
 Let $\Gamma$ be an NEC group of signature $(0; +; [6,6]; \{(5,8,12)\})$ and presentation with generators 
 
 $$x_1,x_2,e,c_0,c_1,c_2,c_3,$$ and relations
 
 $$x_1^6=x_2^6=c_0^2=c_1^2=c_2^2=c_3^2= (c_0c_1)^5=(c_1c_2)^8=(c_2c_3)^{12}=ec_0e^{-1}=x_1x_2e=1$$
 
 Consider the following permutation representation of $\Gamma$
 
 \bigskip
 
 $x_1\mapsto (1,4)(2,3,6)(5)$
 
 \medskip
 
 $x_2\mapsto (1,6,2,5,4,3)$
 
 \medskip
 
 $e\mapsto (3,2,1)(6,5,4))$
 
 $c_0\mapsto (1,2)(3,4)(5)(6)$
 
 \medskip
 
 $c_1\mapsto (1,3)(2,6)(4)(5)$
 
 \medskip
 
 $c_2\mapsto (1,4)(2,6)(3,5)$
 
 \medskip$c_3\mapsto (1,5)(2,3)(4)(6)$

 \medskip
 
 Let $\Lambda$ be the stabiliser of a point. We want to find the signature of $\Lambda$.  We use our algorithm. 
 
 \medskip
 
 Step 1. The induced reflections are  $c_{05}, c_{06}, c_{14}, c_{15}, c_{34}, c_{36}.$
 
 \medskip
 
 Step 2. The links are $c_0\sim c_1$ with link period 5, $c_1\sim c_2$ with link period  with link period 8, $c_2\sim c_3$ with link period 12, $e^{-1}c_0e\sim c_3$ with link period 1.
 
 \medskip

 Step 3. Find the orbits of the dihedral subgroups.
 
 $\langle c_0,c_1\rangle\cong D_5$,    The orbits are $\{1,2,3,4,6\},\{5\}$
 
 \medskip
 
 $\langle c_1,c_2 \rangle\cong D_8$,    The orbits are $\{1,3,4,5\},\{2,6\}$
 
 \medskip
 
 $\langle c_2,c_3 \rangle\cong D_{12}$,    Orbit $\{1,2,3,4,5,6\}$
 
 \medskip
 
 $\langle ec_0e^{-1}, c_3\rangle\cong D_1$,  The  orbits are $\{1,5\},\{2,3\},\{4\}, \{6\}$
 
 \bigskip
 
 Step 4. On the orbit $\{1,2,3,4,6\}$ of $\langle c_0,c_1 \rangle$ we have one fixed point namely 6 of $c_0$ and one fixed point 5 of $c_1$. This gives us a link $c_{06}\sim c_{15}$ with link period 5/5=1.  On the orbit $\{5\}$ we have one fixed point of both $c_0$ and $c_1$ This gives us a link $c_{05}\sim c_{15}$ with link period 5/1=5.
 
 \medskip
 
 On the orbit $\{1,3,4,5\}$ of $\langle \{c_1,c_2\}\rangle$ we have two fixed points namely 4 and 5 of $c_1$, giving the link $c_{14}\sim c_{15}$ . Here the period of the elliptic element is 8 and $c_1c_2$ has period 4 so that the link period is $8/4=2$.
 
 Thus we have the following links, which give us one chain,
 
 $$c_{00}\sim c_{15}\sim c_{14}\sim c_{06}\sim c_{34}\sim c_{36}\sim c_{05}.$$
 
The link periods are 5 between $c_{05}$ and $c_{15}$, 2 between $c_{15}$ and $c_{14}$, 1 between $c_{14}$ and $c_{06}$, 1 between $c_{06}$ and $c_{34}$, 2 between $c_{34}$ and $c_{36}$ and 1 between $c_{36}$ and $c_{05}.$

This gives us one period cycle $\{(2,2,5)\}$.

The orbit $\{2,6\}$ of $\langle {c_1,c_2} \rangle $ contains no fixed points of $c_1$ or $c_2$ and so by step 4 of the algorithm  this gives an elliptic period $8/1=8$  on $\Lambda$. Also, the orbit $\{2,3\}$ of $\langle ec_0e^{-1}, c_3\rangle$ contains no fixed points of $c_3$ and so gives an elliptic period 1 (which we can ignore).
 
 Also, the elliptic elements $x_1,x_2$ induce elliptic periods on $\Lambda$ of periods 2,3, and 6 from the permutation theorem for Fuchsian groups, \cite{S1}.
 
 Thus we have found that $\Lambda$ has elliptic periods [2,3,6,8] and one period cycle $\{(2,2,5)\}$.

For the orientability we draw the augmented Schreier coset graph. We find an edge labelled $c_0$ joining 1 and 2, an edge labelled $c_3$ joining 2 and 3 and an edge labelled $c_1$ joining 3 and 1. This gives us a triangle in the augmented Schreier graph and hence $\Lambda$ has non-orientable quotient space,
We now use the Riemann-Hurwitz formula as described above to find the genus of $\Lambda$. We find that it is 9, and so $\Lambda$ has NEC signature
$$(9:-;[2,3,6,8];\{(2,2,5)\}).$$
 
 \medskip

\end{document}